\documentclass[aps,prl,reprint,superscriptaddress]{revtex4-2}

\usepackage{subcaption}
\usepackage{amsmath,amssymb,mathtools}
\usepackage{graphicx}
\usepackage{bm}
\usepackage{pgf}
\usepackage{tikz}
\usepackage{hyperref}
\usepackage{float}

\graphicspath{{figs/}{./}}

\begin{document}

\title{Exact Conservative Dynamics from Quantized Interaction Rules}

\author{Park Junhu}
\email{nller.park@snu.ac.kr}
\affiliation{Seoul National University}
\affiliation{iTrix Co., Ltd}

\author{Youngsoo Ha}
\affiliation{Seoul National University}

\author{Myungjoo Kang}
\email{mkang@snu.ac.kr}
\affiliation{Seoul National University}
\affiliation{iTrix Co., Ltd}

\begin{abstract}

Conservation laws are conventionally discretized through floating-point flux evaluation, with invariants obtained by cancellation of approximate interface contributions and admissible weak solutions selected by reconstruction and Riemann solvers. Here we introduce an operator-level formulation in which conservative dynamics is realized as an exact discrete interaction rule on a quantized state space. The update is defined by an antisymmetric integer-transfer operator, which enforces conservation exactly at the arithmetic level and eliminates round-off drift from the primitive evolution \cite{highamAccuracyStabilityNumerical2002}. For scalar laws, monotone order-preserving transfers select admissible shock structures within the primitive update, rather than through flux reconstruction. Numerical experiments show that the interaction rule preserves high-frequency transport near the Nyquist limit and maintains sharply localized discontinuities in Burgers dynamics. The same construction extends to multidimensional problems and systems of conservation laws through oriented, vector-valued integer transfers. The results show that exact integer-transfer dynamics can suppress cumulative transport drift while preserving entropy-shock localization in nonlinear conservative evolution.

\end{abstract}

\maketitle

\section{Introduction}

Conservation laws play a central role in physical modeling, expressing invariant quantities such as mass, momentum, and energy through local balance relations. In standard numerical formulations, these laws are realized by floating-point evaluation of interface fluxes in finite-volume or finite-difference frameworks \cite{toroRiemannSolversNumerical2009,leveque_finite_2002,hartenHighResolutionSchemes1997,osherHighResolutionSchemes1984,haSixthorderWeightedEssentially2016,haImprovingAccuracyFifthOrder2021}. Conservation follows from pairwise cancellation of flux contributions, while admissible weak solutions are selected by reconstruction procedures and Riemann solvers.

Structurally, this approach represents conservative dynamics through approximate flux evaluation in finite-precision arithmetic. Consequently, invariants are preserved only up to machine precision, and the primitive update depends on rounding behavior, summation order, and reconstruction choices. Although these effects are usually small, they create a separation between the formal conservation law and its discrete realization.

In this work, we take a different viewpoint: conservative dynamics is realized directly as a discrete interaction process. Instead of computing approximate fluxes, we define the primitive evolution by a local antisymmetric transfer operator acting on a quantized state space \cite{shannonMathematicalTheoryCommunication1948, coverElementsInformationTheory2006}. Conservation is then not imposed by numerical cancellation, but is an exact algebraic invariant of the interaction rule. The resulting method does not approximate a flux-based scheme at the primitive level; it replaces the flux representation itself by an exact antisymmetric transfer operator on a countable state space.

This work instead formulates a primitive discrete conservative dynamics whose coarse-grained behavior reproduces standard continuum conservation laws.

This representation has two immediate consequences. First, conservation is exact at every step and is independent of floating-point arithmetic. Second, weak-solution selection is relocated to the primitive update: entropy-admissible solutions are determined by the monotone, order-preserving structure of the transfer operator itself \cite{crandallRelationsNonexpansiveOrder1980}. In this sense, the role usually attributed to Riemann solvers is not removed, but absorbed into the operator-level interaction rule.

We therefore interpret conservative dynamics not as the result of flux balancing, but as the coarse-grained description of an underlying discrete interaction rule in which both conservation and entropy selection are encoded in the primitive operator. The continuum flux appears only as an emergent representation of this interaction, while the primitive dynamics is governed by exact integer transfer.

\section{Quantized conservative update}

We start from the scalar hyperbolic conservation law \cite{evansPartialDifferentialEquations2010}

\begin{equation}
\partial_t u + \partial_x f(u) = 0, \label{eq:conservation}
\end{equation}

although the construction is independent of spatial dimension and of the number of conserved components. In multiple dimensions, an antisymmetric integer transfer is assigned to each oriented cell face; for systems, the integer state and transfer become componentwise vectors. The global integer invariant then follows by pairwise cancellation over internal faces, with boundary terms treated in the usual conservative manner. We restrict attention here to the scalar case to expose the core mechanism.

FQNM evolves a quantized state
\begin{equation}
q_i^n \in \mathbb{Z}, \qquad u_i^n \simeq \delta q_i^n,
\end{equation}
where $\delta$ denotes the quantization scale.

The conservative update is
\begin{equation}
q_i^{n+1}
=
q_i^n
-
\left(
F_{i+\frac12}^n - F_{i-\frac12}^n
\right),
\label{eq:update}
\end{equation}

where the integer interface transfer is
\begin{equation}
F_{i+\frac12}^n
=
\phi^+(q_i^n) + \phi^-(q_{i+1}^n).
\label{eq:transfer}
\end{equation}

A representative construction is
\begin{equation}
\phi^\pm(q)
=
\operatorname{round}\!\left[
\frac{\Delta t}{\Delta x}
\frac{f^\pm(\delta q)}{\delta}
\right].
\label{eq:split-transfer}
\end{equation}

Because each interface transfer enters adjacent cells with opposite signs, Eq.~\eqref{eq:update} guarantees exact conservation:
\begin{equation}
\sum_i q_i^{n+1} = \sum_i q_i^n.
\end{equation}

This identity holds exactly in integer arithmetic and does not rely on cancellation of floating-point quantities. Thus the primitive integer mass is invariant in time and independent of summation order.

The quantization error satisfies
\begin{equation}
    \label{quantization_error}
\|u - \delta q\|_\infty \le \frac{\delta}{2}.
\end{equation}

The estimate \eqref{quantization_error} shows that the quantized representation converges uniformly to the continuum state as $\delta\to0$ at the level of state encoding. At finite $\delta$, however, the primitive evolution remains an exact integer dynamics rather than a floating-point flux approximation. The numerical tests below therefore focus on finite-resolution consequences of the operator rule: exact conservation, high-frequency transport coherence, and entropy-anchor localization in nonlinear dynamics.

\section{Numerical tests}

We test whether the integer-transfer formulation preserves high-frequency transport and nonlinear shock structure without introducing cumulative drift or excessive smoothing. All simulations use a uniform grid. The comparisons are designed to separate three effects: transport drift in high-frequency advection, viscous smoothing in a finite-difference baseline, and nonlinear reconstruction effects in an inviscid WENO5--RK3 shock-capturing baseline. In the Burgers tests, all diagnostics are evaluated against the same Hopf--Lax entropy reference and at the same fixed entropy anchors.

\subsection{Gaussian high-frequency packet}

\begin{figure}[H]
\centering
\resizebox{0.95\linewidth}{!}{\input{figs/packet_center_drift_kdx.pgf}}
\caption{Packet-center drift for transported high-frequency Gaussian packets. The vertical axis is the absolute displacement of the packet energy center relative to the exact solution, $|x_c(u)-x_c(u_{\mathrm{exact}})|/\Delta x$, measured in grid-cell units. Here $x_c$ denotes the circular centroid of the mean-subtracted packet energy density $(u-\bar u)^2$, and the horizontal axis is the normalized carrier wavenumber $k_0\Delta x$. Across the tested grid resolutions and carrier frequencies, FQNM keeps the packet-center drift at the sub-cell level, typically $10^{-3}$--$10^{-2}$ grid cells, whereas the WENO5--RK3 baseline accumulates substantially larger displacement and reaches multi-cell drift near the Nyquist regime.}
\label{fig:gaussian}
\end{figure}

Figure~\ref{fig:gaussian} evaluates high-frequency transport through packet displacement rather than through a global norm error. This choice isolates numerical drift of the transported structure from amplitude attenuation: the plotted observable measures how far the packet energy center moves away from the exact advected packet, in units of the grid spacing. The separation is pronounced. The quantized transfer dynamics keeps the packet center locked to the exact trajectory to well below one grid cell throughout the tested regime, while the WENO5--RK3 baseline shows systematic positional drift that grows sharply for near-grid-scale packets. Thus the high-frequency test directly supports the central claim that the primitive integer-transfer rule suppresses transport drift rather than merely reducing a global error norm.

\subsection{Burgers shock formation}

\begin{figure*}[t]
\centering
\begin{subfigure}{0.32\textwidth}
\centering
\resizebox{\linewidth}{!}{\input{figs/snapshot_zoom_t0p140.pgf}}
\caption{$t=0.140$}
\end{subfigure}
\begin{subfigure}{0.32\textwidth}
\centering
\resizebox{\linewidth}{!}{\input{figs/snapshot_zoom_t0p160.pgf}}
\caption{$t=0.160$}
\end{subfigure}

\begin{subfigure}{0.32\textwidth}
\centering
\resizebox{\linewidth}{!}{\input{figs/snapshot_zoom_t0p250.pgf}}
\caption{$t=0.250$}
\end{subfigure}

\begin{subfigure}{0.32\textwidth}
\centering
\resizebox{\linewidth}{!}{\input{figs/snapshot_zoom_t1p000.pgf}}
\caption{$t=1.000$}
\end{subfigure}

\caption{Time-resolved Burgers shock-anchor zooms in the fixed window $x\in[0.475,0.525]$. The initial condition is $u_0(x)=\sin(2\pi x/L)$, for which the entropy shock is fixed by symmetry at $(x,u)=(0.5,0)$. Each panel compares the Hopf--Lax entropy solution, the viscous finite-difference baseline, the FQNM integer-transfer dynamics, and the inviscid WENO5--RK3 baseline before and after shock formation. Markers denote the fixed entropy anchors used in the diagnostics: the central shock anchor and the adjacent left/right entropy cells.}
\label{fig:burgers}
\end{figure*}

\subsection{Fixed-anchor slope error}

\begin{figure}[H]
\centering
\resizebox{0.95\linewidth}{!}{\input{figs/slope_error_vs_fixed_anchors.pgf}}
\caption{One-sided slope error at fixed entropy anchors in Burgers dynamics. The left entropy anchor is evaluated using only its left-sided slope, and the right entropy anchor using only its right-sided slope; the shock anchor itself is not differentiated across. Errors are measured relative to the Hopf--Lax entropy solution at the same fixed anchor points. After shock formation, FQNM rapidly suppresses the entropy-anchor slope error, whereas WENO5--RK3 retains an $O(1)$ local slope discrepancy over long times.}
\label{fig:slope-error}
\end{figure}

The Burgers tests use the initial condition $u_0(x)=\sin(2\pi x/L)$, for which the Hopf--Lax solution is constructed from the corresponding primitive.

Figure~\ref{fig:burgers} probes entropy-shock localization directly at the grid scale over time. Instead of relying on a single final-time snapshot, the fixed-window sequence resolves the shock midpoint and the neighboring entropy cells before and after shock formation. The anchor point $(0.5,0)$ is imposed by symmetry for the sinusoidal initial condition and is therefore not re-detected from a numerical gradient; it provides a fixed diagnostic for transport drift and shock displacement.

The separation between the schemes is most visible after the entropy shock forms. The viscous finite-difference baseline spreads the transition over multiple cells. The inviscid WENO5--RK3 baseline captures a sharp discontinuity, but retains a persistent local mismatch near the fixed entropy anchors. In contrast, the quantized transfer dynamics remains aligned with the midpoint structure and the neighboring entropy-cell configuration.

Figure~\ref{fig:slope-error} quantifies this effect by comparing one-sided slopes at the same fixed anchors. The diagnostic deliberately avoids differentiating across the discontinuity: the left anchor is compared using only its left-sided slope and the right anchor using only its right-sided slope. After shock formation, the FQNM slope error rapidly decreases and remains small, whereas the WENO baseline retains an $O(1)$ discrepancy over long times. This supports the interpretation that the quantized transfer rule suppresses cumulative entropy-anchor drift rather than merely producing a visually sharp shock profile.

The remaining finite-grid difference from the sampled Hopf--Lax profile should be interpreted as compatibility between a continuum entropy reference and a quantized conservative update. The Hopf--Lax construction gives the continuum entropy solution, whereas the primitive rule evolves an integer grid state satisfying Eq.~\eqref{eq:update} exactly at each step. The fixed-anchor diagnostics therefore compare the local entropy-shock geometry, not pointwise equality of the continuum and quantized profiles.

\section{Discussion}

The present formulation replaces flux-based reconstruction, at the primitive level, with a discrete antisymmetric transfer operator. Invariance is enforced algebraically by the update itself and is therefore independent of floating-point round-off.

Extensions to multidimensional problems and systems follow from the same operator-level construction, although they are not developed here in detail. The present results suggest that exact integer-transfer dynamics can reproduce conservative continuum behavior while suppressing cumulative transport drift at finite resolution.

\section{Conclusion}

We have introduced an operator-level formulation of conservative dynamics based on exact integer transfer between quantized states. The resulting update conserves exactly in integer arithmetic, avoids round-off drift, and encodes entropy selection through a monotone order-preserving transfer rule. Numerical tests show that this structure preserves high-frequency transport and entropy-shock localization while suppressing cumulative transport drift in nonlinear conservative evolution.

\bibliographystyle{apsrev4-2}
\bibliography{refs}

\end{document}